\documentclass[12pt]{amsart}
\usepackage{amsmath, amssymb,graphicx,caption, subcaption, mathtools, enumitem}
\usepackage[pdftex,hyperindex=true]{hyperref}
\usepackage{url}
\usepackage{etoolbox}
\patchcmd{\section}{\scshape}{\bfseries}{}{}
\patchcmd{\subsection}{\bfseries}{\itshape}{}{}

\makeatletter
\def\@seccntformat#1{%
  \protect\textup{\protect\@secnumfont
    \ifnum\pdfstrcmp{section}{#1}=0 \bfseries\fi
    \ifnum\pdfstrcmp{subsection}{#1}=0 \itshape \fi
    \csname the#1\endcsname
    \protect\@secnumpunct
  }%
}  
\makeatother

\usepackage[utf8]{inputenc}
\usepackage[english]{babel}

\DeclarePairedDelimiter\abs{\lvert}{\rvert}%
\theoremstyle{plain}

\newtheorem{theorem}{Theorem}[section]

\newtheorem*{conjecture}{Conjecture}

\theoremstyle{definition}

\newtheoremstyle{note}
{3pt}
{3pt}
{\itshape}
{}
{\itshape}
{:}
{.5em}
{}

\theoremstyle{note}
\newtheorem{remark}[theorem]{Remark}

\theoremstyle{plain} 
\newcommand{\thistheoremname}{}
\newtheorem*{genericthm}{\thistheoremname}

\numberwithin{equation}{section}

\newcommand{\acknowledge}{\subsection*{Acknowledgements}}

\def\Dbar{\leavevmode\lower.6ex\hbox to 0pt{\hskip-.23ex \accent"16\hss}D}

\begin{document}

\title{Weights, Weyl-equivariant maps and a rank conjecture}
\author{J. Malkoun}
\address{Department of Mathematics and Statistics\\
Faculty of Natural and Applied Sciences\\
Notre Dame University-Louaize, Zouk Mosbeh\\
P.O.Box: 72, Zouk Mikael,
Lebanon}
\email{joseph.malkoun@ndu.edu.lb}

\date{\today}
\maketitle

\vspace{0.3cm}
\centerline{\em{In memory of Sir Michael Atiyah}}
\vspace{0.3cm}

\begin{abstract} In this note, given a pair $(\mathfrak{g}, \lambda)$, where $\mathfrak{g}$ is a complex semisimple Lie algebra and $\lambda \in \mathfrak{h}^*$ is a 
dominant integral weight of $\mathfrak{g}$, where $\mathfrak{h} \subset \mathfrak{g}$ is the real span of the coroots inside a fixed Cartan subalgebra, we associate 
an $SU(2)$ and Weyl equivariant smooth map $f: X \to (P^m(\mathbb{C}))^n$, where $X \subset \mathfrak{h} \otimes \mathbb{R}^3$ is the configuration 
space of regular triples in $\mathfrak{h}$, and $m$, $n$ depend on the initial data $(\mathfrak{g}, \lambda)$.

We conjecture that, for any $\mathbf{x} \in X$, the rank of $f(\mathbf{x})$ is at least the rank of a collinear configuration in $X$ (collinear when viewed as an ordered $r$-tuple of points in $\mathbb{R}^3$, with $r$ being 
the rank of $\mathfrak{g}$). A stronger conjecture is also made using the singular values of a matrix representing $f(\mathbf{x})$.

This work is a generalization of the Atiyah-Sutcliffe problem to a Lie-theoretic setting.
\end{abstract}

\maketitle

\section{Introduction} \label{intro}

While its origin lies in Physics, more specifically in the work of Berry and Robbins \cite{BR1997} on a geometric explanation of the spin-statistics theorem, the 
Atiyah-Sutcliffe problem on configurations of points is a geometric problem. Consider
\[ C_n(\mathbb{R}^3) = \{ \mathbf{x} = (\mathbf{x}_1,\ldots,\mathbf{x}_n) \in (\mathbb{R}^3)^n; \mathbf{x}_a \neq \mathbf{x}_b \text{ for all $a \neq b$} \} \]
where $1 \leq a,b \leq n$. Also consider the flag manifold $U(n)/T^n$. We note that the symmetric group $\Sigma_n$ on $n$ elements acts on 
$C_n(\mathbb{R}^3)$ by permuting the $n$ points $\mathbf{x}_a \in \mathbb{R}^3$, for $1 \leq a \leq n$. Moreover, $\Sigma_n$ also acts on 
the flag manifold $U(n)/T^n$, thought of as a left coset space, by permuting the $n$ columns of a matrix $g \in U(n)$ representing a left 
coset $gT^n \in U(n)/T^n$.

The Berry-Robbins problem asks whether there exists, for each $n \geq 2$, a continuous map
\[ f_n: C_n(\mathbb{R}^3) \to U(n)/T^n \]
which is $\Sigma_n$ equivariant.

The Berry-Robbins problem was first solved positively by M.F. Atiyah in \cite{Ati-2000}. However, the solution there had some unsatisfactory features. In that same article, 
and a subsequent article \cite{Ati-2001}, other candidate maps $f_n$ were presented, with more satisfactory features (for instance, these candidate maps are 
\emph{smooth}), but they would be genuine solutions only provided some linear independence conjecture holds. Later in \cite{Ati-Sut-2002}, Sir Michael Atiyah and 
Paul Sutcliffe found strong numerical evidence for linear independence, as well as for a stronger conjecture, which says that $\abs{D(\mathbf{x})} \geq 1$ for any 
$\mathbf{x} \in C_n(\mathbb{R}^3)$, where $D: C_n(\mathbb{R}^3) \to \mathbb{C}$ is a (smooth) normalized determinant function whose non-vanishing is equivalent to 
the linear independence conjecture. These are the Atiyah-Sutcliffe conjectures $1$ and $2$ respectively. The authors of \cite{Ati-Sut-2002} also formulated a conjecture $3$ 
which implies conjecture $2$, but we will not explain it here.

In \cite{Ati-Bie}, M.F. Atiyah and R. Bielawski solved a Lie-theoretic generalization of the Berry-Robbins problem using Nahm's equations. However, their solution was not explicit. 
In his Edinburgh Lectures on Geometry, Analysis and Physics \cite{Edin_lectures}, Sir Michael asked whether there exists a Lie-theoretic generalization of the Atiyah-Sutcliffe problem. 
In this article, we provide a positive answer to his question, and generalize the Atiyah-Sutcliffe problem to a Lie-theoretic setting.

\section{A Lie-theoretic generalization of the Atiyah-Sutcliffe problem} \label{generalization}

Let $\mathfrak{g}$ be a complex semisimple Lie algebra. Denote its Killing form by $(-,-)$.  Let $\mathfrak{h} \subset \mathfrak{g}$ be 
the real span of the coroots of $\mathfrak{g}$ inside a fixed Cartan subalgebra (the latter being thus $\mathfrak{h} \otimes \mathbb{C}$). Let $R \subset \mathfrak{h}^*$ be the set of 
all roots of $\mathfrak{g}$ with respect to $\mathfrak{h} \otimes \mathbb{C}$, denote by $R^+ \subset \mathfrak{h}^*$ a choice of positive roots, 
and by $\Phi \subset \mathfrak{h}^*$ the corresponding set of simple roots. Strictly speaking, a root $\alpha$ of $\mathfrak{g}$ is an element of 
$(\mathfrak{h} \otimes \mathbb{C})^*$ but, since each root $\alpha$ is real-valued on $\mathfrak{h}$, we consider each $\alpha$ as an element of $\mathfrak{h}^*$. 
It is known that the Killing form (up to a sign) on $\mathfrak{g}$ restricts to an inner product on $\mathfrak{h}$.

Denote by $W$ the Weyl group of $(\mathfrak{g}, \mathfrak{h})$. Thus $W$ is the group generated by reflections 
$r_{\alpha}$ in $\mathfrak{h}^*$ with respect to the hyperplane $\alpha^{\perp}$, as $\alpha$ varies in the set of roots $R$. We let 
$\lambda \in \mathfrak{h}^*$ be a dominant integral weight of $\mathfrak{g} \otimes \mathbb{C}$. What the latter condition amounts to is that
\[ 2 \frac{(\alpha, \lambda)}{(\alpha,\alpha)} \]
is a nonnegative integer for any positive root $\alpha \in R^+$.

Denote by $X$ the following configuration space
\[ X = \{ \mathbf{x} \in \mathfrak{h} \otimes \mathbb{R}^3; (\alpha \otimes 1)(\mathbf{x}) \neq \mathbf{0} \text{ for any $\alpha \in R^+$} \} \]
where $\alpha \otimes 1: \mathfrak{h} \otimes \mathbb{R}^3 \to \mathbb{R}^3$ is the linear map obtained by tensoring $\alpha$ with the identity map 
on $\mathbb{R}^3$.

From now on, we identify $\mathfrak{h}$ with $\mathfrak{h}^*$ via the (restricted) Killing form $(-,-)$, so that the Weyl group $W$ acts naturally on $\mathfrak{h}$. If we tensor this action with the trivial 
action of $W$ on $\mathbb{R}^3$, we obtain an action of $W$ on $\mathfrak{h} \otimes \mathbb{R}^3$, which preserves $X$. On the other hand, $SU(2)$ acts on 
$\mathfrak{h} \otimes \mathbb{R}^3$ via the tensor product of the trivial action on $\mathfrak{h}$ and its natural action on $\mathbb{R}^3$ via its adjoint action, i.e. 
the $2$-to-$1$ group homomorphism from $SU(2)$ onto $SO(3)$. This action of $SU(2)$ preserves $X$.

We let
\[ m = \sum_{\alpha \in R^+} 2 \frac{(\alpha, \lambda)}{(\alpha,\alpha)} \]
and
\[ n = [W:W_\lambda] ,\]
where $W_{\lambda}$ is the stabilizer of $\lambda$ in $W$.

Given the initial data $(\mathfrak{g}, \mathfrak{h}, R^+, \lambda)$ as above, we will construct a smooth map $f: X \to (P^m(\mathbb{C}))^n$. Let $\mathbf{x} \in X$. 
For any root $\alpha \in R$, we define $v_{\alpha} \in S^2 \subset \mathbb{R}^3$ as the normalization (with respect to the Euclidean inner product on $\mathbb{R}^3$) of
\[ (\alpha \otimes 1)(\mathbf{x}) \in \mathbb{R}^3 \setminus \{ \mathbf{0} \} \]

The Hopf map $h: S^3 \to S^2$ can be defined by
\[ h(u,v) = (2 u \bar{v}, \abs{u}^2 - \abs{v}^2) \]
where $S^3 \subset \mathbb{C}^2$ is the set of all $(u,v) \in \mathbb{C}^2$ such that $\abs{u}^2 + \abs{v}^2 = 1$, and 
$S^2 \subset \mathbb{C} \times \mathbb{R}$ is defined as the set of all $(\zeta,z) \in \mathbb{C} \times \mathbb{R}$ such that $\abs{\zeta}^2 + z^2 = 1$.

For every root $\alpha \in R$, we choose a Hopf lift $(u_{\alpha}, v_{\alpha}) \in S^3$. Such a Hopf lift is unique up to a global factor in $U(1)$. 
We then form the complex polynomial
\[ p_{\alpha}(t) = u_{\alpha} t - v_{\alpha} .\]

The elements of $W/W_{\lambda}$ are in natural one-to-one correspondence with the Weyl orbit $W.\lambda$ of $\lambda \in \mathfrak{h}^*$. Let us say that
\[ W.\lambda = \{ \lambda_1, \ldots, \lambda_n \} .\]
Choose $g_1,\ldots, g_n \in W$ so that $g_k(\lambda) = \lambda_k$, for $1 \leq k \leq n$.

For any $1 \leq k \leq n$, let
\[ p_k(t) = \prod_{\alpha \in R^+} ( p_{g_k.\alpha}(t) )^{m_{\alpha}} \]
where
\[ m_{\alpha} = 2 \frac{(\alpha, \lambda)}{(\alpha,\alpha)} .\]
The latter is a nonnegative integer since $\lambda$ is a dominant integral weight of $\mathfrak{g} \otimes \mathbb{C}$ and $\alpha$ is a positive root. We remark that 
the definition of $p_k$ does not depend on the choice of representative $g_k$ in its left coset $g_k W_{\lambda}$, since another choice, say $g_k w$, where 
$w \in W_{\lambda}$, will only permute the factors of $p_k$, since the map which maps $\alpha$ to $w.\alpha$ permutes the set of positive roots 
$\alpha$ satisfying $(\alpha, \lambda) > 0$.

We define the map
\[ F: X \to ( \mathbb{C}^{m+1} \setminus \{ \mathbf{0} \} )^n / T^n \]
which maps each $\mathbf{x} \in X$ to the $n$-tuple of polynomials $(p_k(t))$, for $1 \leq k \leq n$, where each 
polynomial is actually only defined up to multiplication by an element of $U(1)$ (due to the $U(1)$ ambiguity of each Hopf lift). 
Finally, the map
\[ f: X \to ( P^m(\mathbb{C}) )^n \]
is obtained by following the map $F$ with the natural projection
\[ ( \mathbb{C}^{m+1} \setminus \{ \mathbf{0} \} )^n / T^n \to ( P^m(\mathbb{C}) )^n .\]

We note that an element $\mathbf{x} \in \mathfrak{h} \otimes \mathbb{R}^3$ can be thought of as an $r$-tuple of points in $\mathbb{R}^3$, where 
$r = \operatorname{rank}(\mathfrak{g}) = \operatorname{dim}(\mathfrak{h})$. We then define the class of collinear configurations, as the set of all $\mathbf{x} \in X$ 
consisting of $r$ collinear points in $\mathbb{R}^3$. Given $\mathbf{x} \in X$, $f(\mathbf{x})$ 
can be represented by an $m+1$-by-$n$ complex matrix. Such a choice is not unique, as one can independently scale each column of such a matrix. It turns out 
that the rank of $f(\mathbf{x})$, that is to say the rank of an $m+1$-by-$n$ complex matrix representing $f(\mathbf{x})$, is the same for all collinear configurations 
$\mathbf{x}$. We denote this rank by $r_{col}$.

Our first conjecture can now be phrased. 
\begin{conjecture}[Generalized Conjecture $1$] Given any $\mathbf{x} \in X$, 
\[ \operatorname{rank}(f(\mathbf{x})) \geq r_{col}. \] 
\end{conjecture}

Our second conjecture can be phrased using the singular 
values of a matrix representing $F(\mathbf{x})$, and is a quantitative refinement of our first conjecture. More specifically, define the $m+1$-by-$m+1$ matrix $g$ by
\[ g = \operatorname{diag}\left(\binom{m}{0}^{-1}, \binom{m}{1}^{-1}, \ldots, \binom{m}{m-1}^{-1}, \binom{m}{m}^{-1}\right) \]
and $\Delta: X \to \mathbb{R}$ by
\[ \Delta(\mathbf{x}) = s_{r_{col}} ((\sqrt{g})^{-1} \operatorname{Sing}( \sqrt{g} F(\mathbf{x}) )) \]
where $s_j$ denotes the $j$-th elementary symmetric polynomial of the diagonal entries of the matrix argument (which has $0$s off the diagonal), 
and $\operatorname{Sing}$ denotes the middle matrix in 
the singular value decomposition, namely the matrix containing the singular values (and possibly zero(s)) as diagonal entries, with 
multiplicity taken into account. We note that the set of singular values does not depend on the choice of $m+1$-by-$n$ complex matrix representing 
$F(\mathbf{x})$, since another such matrix is obtained from the first by multiplying by a diagonal unitary matrix from the right.

The matrix $\sqrt{g}$ and its inverse were used in the previous definition in order to make $\Delta$ $SU(2)$-invariant.

We remark that $\Delta$ is clearly non-negative. Moreover, its non-vanishing on $X$ is equivalent to our Generalized Conjecture $1$. We can now formulate 
our Generalized Conjecture $2$.
\begin{conjecture}[Generalized Conjecture $2$] Given any $\mathbf{x} \in X$,
\[ \Delta(\mathbf{x}) \geq \Delta(\mathbf{x}_{col}) .\]
\end{conjecture}

These two conjectures are generalizations of the Atiyah-Sutcliffe conjectures $1$ and $2$ to a Lie theoretic setting. Indeed, if $\mathfrak{g} = \mathfrak{sl}(n)$, 
and
\begin{align*} \lambda &= e^1 + \cdots + e^{n-1} - \frac{n-1}{n}(e^1 + \cdots + e^n) \\
         &= \frac{1}{n}( e^1 + \cdots + e^{n-1} - (n-1) e^n ) 
\end{align*}
where $e_k$ represents the diagonal $n$-by-$n$ matrix having $1$ at the $(k,k)$-entry, and $0$ everywhere else, $\mathfrak{h}$ is given by
\[ \mathfrak{h} = \operatorname{span}_{\mathbb{R}}( e_k - e_{k+1}; \text{ for $1 \leq k \leq n-1$}) \]
and $(e^1, \ldots, e^n)$ represents the dual basis of the basis $(e_1, \ldots, e_n)$ of the space of diagonal matrices. Then 
for such choices, our map $f$ specializes to the Atiyah-Sutcliffe 
map, and our two conjectures specialize to the Atiyah-Sutcliffe conjectures $1$ and $2$. This is thus a generalization of the Atiyah-Sutcliffe problem to 
a Lie-theoretic setting.

\begin{remark} The author's $Sp(m)$ version of the Atiyah-Sutcliffe problem in \cite{Malkoun2014} is also a special case of our Lie-theoretic generalization, and can be obtained 
by a suitable choice of weight $\lambda$ for the Lie algebra $\mathfrak{sp}(2m, \mathbb{C})$.
\end{remark}

\begin{remark} While Atiyah and Bielawski in \cite{Ati-Bie} have found a Lie-theoretic solution to the so-called Berry-Robbins problem, it is not clear, as 
of now, how their non-elementary solution, which uses Nahm's equations, is related to the more elementary Atiyah-Sutcliffe problem. It would be very interesting to try 
and relate the two approaches, if possible.
\end{remark}

\section{Numerical Evidence}

The author did some numerical testing of the Generalized Conjecture $2$ for $\mathfrak{g} = \mathfrak{sl}(4)$ for $26$ different weights $\lambda$. For each such a 
weight, the computer generated $1000$ configurations pseudo-randomly for which it calculated their $\Delta$'s, and then calculated the minimal $\Delta$ 
among these $1000$ configurations. Such a sample-minimum $\Delta$ is shown in the table below, for these $26$ weights, as well as the corresponding 
$\Delta(\mathbf{x}_{col})$. We can see that in all these cases, the former is greater or equal to the latter, thus supporting our Generalized Conjecture $2$.

\begin{table}[h!]
  \begin{center}
    \caption{Numerical Tests for $\mathfrak{sl}(4)$}
    \label{tab:table1}
    \begin{tabular}{l | l | l} 
      \textbf{weight} & \textbf{sample-min. $\Delta$} & \textbf{collinear $\Delta$}\\
      \hline
      $[6, 4, 2, 0]$ & 13697132122.474367 & 1086.1160159029314 \\
 $[5, 3, 1, 0]$ & 519458.6316778218 & 56601.99402847759 \\
 $[4, 2, 0, 0]$ & 61189.2491373496 & 45.25483399593905 \\
 $[5, 3, 2, 0]$ & 2340879.6536430004 & 34796.689497708816 \\
 $[4, 2, 1, 0]$ & 202.89393151348898 & 23.99999999999951 \\
 $[3, 1, 0, 0]$ & 1050.3074380238 & 107.33126291998899 \\
 $[4, 2, 2, 0]$ & 12577.200441098057 & 48.00000000000035 \\
 $[3, 1, 1, 0]$ & 5.1901705590915075 & 1.9999999999999858 \\
 $[2, 0, 0, 0]$ & 4.1732831617086825 & 2.828427124746185 \\
 $[5, 4, 2, 0]$ & 1292354.6342256288 & 56601.99402847781 \\
 $[4, 3, 1, 0]$ & 1003499.9746244224 & 28676.856731517517 \\
 $[3, 2, 0, 0]$ & 654.9546591874636 & 99.49874371066126 \\
 $[4, 3, 2, 0]$ & 109.5863462909002 & 23.999999999999552 \\
 $[3, 2, 1, 0]$ & 63.03506766070143 & 33.94112549695389 \\
 $[2, 1, 0, 0]$ & 4.184180177236805 & 3.9999999999999756 \\
 $[3, 2, 2, 0]$ & 3.246186555208909 & 1.9999999999999867 \\
 $[2, 1, 1, 0]$ & 34.135503300239776 & 26.83281572999739 \\
 $[1, 0, 0, 0]$ & 1.0308769367806199 & 0.999999999999999 \\
 $[4, 4, 2, 0]$ & 16337.905038101348 & 45.25483399593978 \\
 $[3, 3, 1, 0]$ & 535.2649778926265 & 99.49874371066095 \\
 $[2, 2, 0, 0]$ & 23.755650809146687 & 5.656854249492361 \\
 $[3, 3, 2, 0]$ & 414.7813428343056 & 107.33126291998886 \\
 $[2, 2, 1, 0]$ & 5.212788207618033 & 3.999999999999977 \\
 $[1, 1, 0, 0]$ & 1.445559164471774 & 1.4142135623730918 \\
 $[2, 2, 2, 0]$ & 3.7560678226481286 & 2.8284271247461827 \\
 $[1, 1, 1, 0]$ & 1.0414407239858756 & 0.999999999999999
    \end{tabular}
  \end{center}
\end{table}

We remark that the notation $[3,3,2,0]$ corresponds to the orthogonal projection of $3e^1 + 3e^2 + 2e^3$ onto the orthogonal complement of 
$e^1 + e^2 + e^3 + e^4$ corresponding to the condition of being tracefree. The simulation above took about $30$ minutes on a Macbook Pro $2015$. We 
wish to run more numerical simulations in the future.

\acknowledge{I dedicate this work to Sir Michael Atiyah who came up with the original problem, as well as the question which motivated this work. The author thanks 
Ben Webster and James Humphreys for their comments on the Mathematics StackExchange website and by email. Any possible mistake in this work is however 
only the author's responsibility.}

\smallskip

\vspace{5mm}

\def\Dbar{\leavevmode\lower.6ex\hbox to 0pt{\hskip-.23ex \accent"16\hss}D}
\providecommand{\bysame}{\leavevmode\hbox to3em{\hrulefill}\thinspace}
\providecommand{\MR}{\relax\ifhmode\unskip\space\fi MR }
\providecommand{\MRhref}[2]{%
  \href{http://www.ams.org/mathscinet-getitem?mr=#1}{#2}
}
\providecommand{\href}[2]{#2}

\end{document}